\documentclass[20pt,a4paper]{article}
\usepackage{color,amsmath,amssymb, lscape,graphicx,ifthen,pdfpages}
\DeclareSymbolFont{boldsymbols}{OMS}{cmsy}{b}{n} 
\DeclareSymbolFontAlphabet{\mathbfcal}{boldsymbols} 
\usepackage[colorlinks=true,
linkcolor=webgreen,
filecolor=webbrown,
citecolor=webgreen]{hyperref}

\definecolor{bittersweet}{rgb}{1.0, 0.44, 0.37}
\definecolor{webgreen}{rgb}{0,.5,0}\definecolor{webbrown}{rgb}{.6,0,0}
\definecolor{antiquebrass}{rgb}{0.8, 0.58, 0.46}
\definecolor{antiquefuchsia}{rgb}{0.57, 0.36, 0.51}
\definecolor{apricot}{rgb}{0.98, 0.81, 0.69}
\definecolor{aquamarine}{rgb}{0.5, 1.0, 0.83}
\definecolor{atomictangerine}{rgb}{1.0, 0.6, 0.4}
\definecolor{auburn}{rgb}{0.43, 0.21, 0.1}
\definecolor{azure(colorwheel)}{rgb}{0.0, 0.5, 1.0}
\definecolor{battleshipgrey}{rgb}{0.52, 0.52, 0.51}
\definecolor{beaublue}{rgb}{0.74, 0.83, 0.9}
\definecolor{bleudefrance}{rgb}{0.19, 0.55, 0.91}

\newcommand{\seqnum}[1]{\href{http://oeis.org/#1}{\underline{#1}}}

\def\dstyle#1{$\displaystyle #1 $}
\def\noin{\noindent}
\def\pn{\par\noindent}

\def\psn{\par\smallskip\noindent}
\def\pbn{\par\bigskip\noindent}
\def\Beq{\begin{equation}}
\def\Eeq{\end{equation}}
\def\Beqarray{\begin{eqnarray}}
\def\Eeqarray{\end{eqnarray}}
\def\Eq#1{eq.\,$(#1)$}
\def\sspgeq{\,\geq} 
\def\sspleq{\, \leq \,}
\def\sspkl{\, < \,}

\def\sspeq{\, =\,}
\def\speq{\ =\ }
\def\sspdef{\, :=\,}
\def\spdef{\ :=\ }
\def\sspfed{\, =:\,}
\def\sspin{\, \in \,}

\def\spp{\ +\ }
\def\spm{\ -\ }
\def\sspp{\, +\ }
\def\sspm{\, -\ }

\def\sspto{\,\to\,}

\def\sspneq{\, \neq \,}

\def\binomial#1#2{{#1} \choose {#2}}
\def\ogf{{\it o.g.f.\ }}
\def\ogfs{{\it o.g.f.}s\ }

\def\egf{{\it e.g.f.\ }}
\def\egfs{{\it e.g.f.}s\ }

\def\lgf{{\it l.g.f.}\ }

\def\fps{{\it f.p.s.}\ } 
\def\fpsb{({\it f.p.s.})\ } 

\def\ie{{\it i.e.},\, }
\def\eg{{\it e.g.},\, }
\def\Eg{{\it E.g.},\, }

\def\sspfed{\, =:\, }

\def\Simn{\,{\lower2pt\hbox{$\buildrel {\lower3pt\hbox{$n$}} \over \sim$}}\,}
\def\simn1#1{\,{\lower2pt\hbox{$\buildrel {\lower3pt\hbox{$#1$}} \over \sim$}}\,}

\def\Chi{\raise2.5pt\hbox{$\chi$}} 

\normalbaselines
\begin{document}
\bibliographystyle{unsrt}
\rightline{Karlsruhe} \par\smallskip\noindent
\rightline{August 04 2017}
\vbox {\vspace{6mm}}
\begin{center}
{\Large {\bf On Generating functions of Diagonals Sequences of Sheffer and Riordan Number Triangles}}\\ [9mm]
Wolfdieter L a n g \footnote{ 
\href{mailto:wolfdieter.lang@partner.kit.edu}{\tt wolfdieter.lang@partner.kit.edu},\quad 
\url{http://www.itp.kit.edu/~wl}
                               } \\[3mm]
\end{center}
\vspace{2mm}
\begin{abstract}
\par\smallskip\noindent
The exponential generating function of ordinary generating functions of diagonal sequences of general {\sl Sheffer} triangles is computed by an application of {\sl Lagrange}'s theorem. For the special {\sl Jabotinsky} type this is already known. An analogous computation for general {\sl Riordan} number triangles leads to a formula for the logarithmic generating function of the ordinary generating functions of the product of the entries of the diagonal sequence of {\sl Pascal}'s triangle and those of the {\sl Riordan} triangle. For some examples these ordinary generating functions yield in both cases coefficient triangles of certain numerator polynomials. 
\end{abstract}
\section{Introduction and Summary}
\hskip 1cm The study of the diagonal sequences of {\sl Sheffer} number triangles (exponential, also known as binomial, lower triangular convolution matrices) is interesting. The name exponential {\sl Riordan} arrays is sometimes used for these triangles. The {\sl Sheffer} structure immediately leads to the exponential generating functions ({\it  e.g.f.}s) of the column sequences. It is more difficult to obtain information about these functions for diagonal sequences. {\sl Bala} \cite{Bala} has shown, following {\sl Drake} \cite{Drake},  for a special type of {\sl Sheffer} triangles, called {Jabotinsky} triangles by {\sl Knuth} \cite{Knuth}, that the \egf of the ordinary generating functions (\ogfs) of the diagonal sequences can be computed from {\sl Lagrange}'s inversion theorem. We present in the first part the result for general {\sl Sheffer} triangles and give some examples. They lead to other number triangles providing the coefficients of the numerator polynomials of the \ogfs of the diagonal sequences. In the second part the same analysis is done for general {\sl Riordan} number triangles (ordinary lower triangular convolution matrices). However, one does not obtain information about the diagonal sequences themselves but on certain products of the diagonal entries with other numbers. We will give the result for the logarithmic generating function of the \ogfs of the sequences of the product of the entries of the diagonals of the {\sl Riordan} and the {\sl Pascal} triangle. (The {\sl Pascal} triangle is a special {\sl Riordan} triangle, and also a special {\sl Sheffer} triangle). Also in this case special examples lead to coefficient triangles for the numerator polynomials of these \ogfs.\psn
For {\sl Sheffer} and {\sl Riordan} triangles see \cite{Roman}, \cite{Shapiro} and the {\sl W. Lang} link \cite{WLang} in {\it OEIS} \cite {OEIS} \seqnum{A006232} (henceforth we will omit the {\it OEIS} reference for A-numbers). There also references can be found.\psn
Proofs for not obvious or not standard {\sl Sheffer} or {\sl Riordan} statements  will be given in section $2$.
\pbn        
{\bf Part A: Sheffer triangles and their diagonals} 
\pbn 
\hskip 1cm A {\sl Sheffer} triangle $S$ (an infinite dimensional lower triangular exponential convolution matrix; for practical purpose a $N\times N$ matrix) is denoted by $S\sspeq (g,\, f)$ with \egf \dstyle{g(s)\sspeq \sum_{k=0}^{\infty}\, g_n\,\frac{s^n}{n!}}, where $g(0) \sspeq g_0 \speq 1$ ({\it w.l.o.g.}), and  $f(s) \sspeq s\,{\hat f}(s)$ with \egf \dstyle{{\hat f}(s)\sspeq \sum_{k=0}^{\infty}\, {\hat f}_n\,\frac{s^n}{n!}}, where ${\hat f}(0)\sspeq {\hat f}_0\sspneq 0$. The column sequence $SCol(m)\sspeq \{S(n, m)\}_{n=0}^{\infty}$ (with $m$ leading zeros) has \egf \dstyle{ESCol(s, m)\sspeq \sum_{n = m}^{\infty}\, S(n, m)\,\frac{s^n}{n!}}, for $m\sspin \mathbb N_0\sspdef \{0,\,1,\,...\}$, given by
\Beq
ESCol(s, m)\sspeq g(s)\,\frac{f(s)^m}{m!}\sspeq g(s)\,\frac{s^m\,{\hat f}(s)^m}{m!}\, .
\Eeq
In this paper formal power series \fpsb are considered, and therefore no convergence issues are treated.
\psn
The (ordinary, not exponential) row polynomials (called {\sl Sheffer} polynomials) are  $PS(n,\, x)\sspeq \sum_{m=0}^n\,S(n,\,m)\,x^n$. They have \egf \dstyle{EPS(s,x)\sspeq \sum_{n=0}^{\infty} PS(n,\, x)\,\frac{s^n}{n!}} given by
\Beq
EPS(s,x)\sspeq g(s)\,e^{x\,f(s)}\,, 
\Eeq
which is also called the \egf of the $S$ triangle.\psn
The important exponential convolution property of {\sl Sheffer} polynomials, implied by \Eq{2}, is
\Beq
PS(n,x\sspp y)\sspeq \sum_{k=0}^n\,{\binomial{n}{k}}\,P(k, x)\, PS(n-k,\,y)\sspeq \sum_{k=0}^n\, {\binomial{n}{k}}\,PS(k, x)\, P(n-k,\,y)\,,
\Eeq
where $P$ are the special {\sl Sheffer} polynomials $P\sspeq (1,\, f)$, called associated polynomials to $S\sspeq (g,\, f)$. (See Roman \cite{Roman} for {\sl Sheffer} sequences of polynomials. The notation there differs from the present one. See the above mentioned {\sl W. Lang} link for the relation between them.)\psn
The diagonal sequences are labeled by $d\sspin \mathbb N_0$, with $d\sspeq 0$ for the main diagonal. Their entries are
\Beq
DS(d,\,m)\sspeq S(d+m,\,m),\ \ {\rm for}\ \ m\sspin \mathbb N_0\,.
\Eeq
Their \ogf is
\Beq
GDS(d,\,t)\sspeq \sum_{m=0}^{\infty}\, DS(d,\,m)\,t^m
\Eeq
( the use of $t$ instead of $x$ is motivated by the later appearance of the parameter $t$), and the \egf of $\{GDS(d,\,t)\}_{d=0}^{\infty}$ is taken as 
\Beq
EGDS(y,\,t)\sspdef \sum_{d=0}^{\infty}\, GDS(d,\,t)\,\frac{y^{d+1}}{(d+1)!}\, .
\Eeq
(The unconventional powers for this \egf and the use of $y$ instead of $s$ will become clear later). \psn
To derive a formula for this \egf $EGDS(y,\,t)$  of \ogfs of diagonal sequences we need {\sl Lagrange}'s theorem and an application.\psn
{\bf Lemma: Lagrange theorem and inversion} \cite{Fichtenholz}, p. 523, \Eq{29},  \cite{WhittakerWatson}, p. 133.\psn
{\bf a)} For ${\widetilde H}(x) = H(y(x))$ with implicit $y\sspeq y(x)\sspeq a\sspp x\,\varphi(y)$  (here as \fps) one has
\Beq
{\widetilde H}(x)\sspeq H(a)\sspp \sum_{n=1}^{\infty}\,\frac{x^n}{n!} \frac{d^{n-1}\ }{da^{n-1}}\left[\varphi^n(a)\,H^{\prime}(a)\right]\, .
\Eeq
{\bf b)} With $a\sspeq 0$, $y\sspeq y(x) \sspeq x\, \psi(x)$, and the compositional inverse $x\sspeq y^{[-1]} \sspeq x(y)$ it follows that
\Beqarray
{\tilde H}(y)\sspeq H(x(y)) &\sspeq& H(0)\sspp \sum_{n=1}^{\infty}\, \frac{y^n}{n!}\,\left.\frac{d^{n-1}\ }{da^{n-1}}\left[\left(\frac{1}{\psi(a)}\right)^n\,H^{\prime}(a)\right]\right|_{a=0} \nonumber \\
&\sspeq&  H(0) \sspp \sum_{n=1}^{\infty}\, \frac{y^n}{n!}\,(n-1)!\,[a^{n-1)}]\,\left[\left(\frac{1}{\psi(a)}\right)^n\,H^{\prime}(a)\right]\,
\Eeqarray
where $[a^n]\,h(a)$ picks the coefficient of $a^n$ of a \fps $h\sspeq h(a)$. 
Applying this {\it Lemma}, part b), introducing a parameter $t$,  to $y\sspeq y(t;x)\sspeq x\,\psi(t;x)\sspeq x\,(1\sspm t\, {\hat f}(x))\sspeq x\sspm t\, f(x)$ with the {\sl Sheffer} function $f$, and taking $H(x)\sspeq \int\,dx\,g(x)$, with the {\sl Sheffer} function $g$, we obtain with the compositional inverse $x\sspeq x(t;y)$ of $y\sspeq y(t;x)$
\pbn
{\bf Proposition 1:}
\Beq
EGDS(y,t) \sspeq H(x(t;y))\sspm H(0)\sspeq \left.\left[\int\,dx\,g(x)\right]\right|_{x\sspeq x(t;y)}\sspm \left.\left [\int\,dx\,g(x)\right]\right|_{x\sspeq 0}\ .
\Eeq
As the last equation shows one has first to compute $x\sspeq x(t;y)$, the compositional ({\sl Lagrange}) inverse of $y\sspeq y(t;x)$. This is the case $H(x) =x$ in the {\it Lemma}, part b, with the chosen $\psi\sspeq \psi(t;x)\sspeq 1-t\,{\hat f}(x)$. This belongs to the associated {\sl Sheffer} case $J\sspeq (1,\,f)$ (the {\sl Jabotinsky} type \cite{Knuth}, here called $J$ instead of $S$). This yields the following corollary which has been treated already by {\sl Bala} \cite{Bala}.\psn
{\bf Corollary 1: Jabotinsky case}\psn
\Beq
EGDJ(y,t) \sspeq x(t;y)\, .
\Eeq
This means that for $J\sspeq (1,\, f)$ the \egf of the \ogfs of the diagonal sequences is just the compositional inverse of $y\sspeq y(t;x)\sspeq x\sspm t\,f(x)$.\pbn
{\bf Examples}
\pbn
{\bf 1)} (\cite{Drake}, Example $1.10.1$, and \cite{Bala}, Example 2) $J\sspeq (1,\, e^s\sspm 1)$, the {\sl Stirling} triangle of the second kind, given in \seqnum{A048993}.
The {\sl Lagrange} inverse  $x\sspeq x(t;y)$ of \dstyle{y\sspeq x\,\left (1\sspm t\,\frac{e^x\sspm 1}{x}\right)} turns out to be (for Maple \cite{Maple} one uses the expansion up to some power to avoid error messages from $x\sspto 0$)
\Beq
 x(t;y)\sspeq \frac{1}{1\sspm t}\,y \sspp \frac{t}{(1\sspm t)^3}\,\frac{y^2}{2!} \sspp \frac{t\,(1\sspp 2\,t)}{(1\sspm t)^5}\,\frac{y^3}{3!} \sspp \frac{t\,(1 \sspp 8\,t \sspp 6\,t^2)}{(1\sspm t)^7}\,\frac{y^4}{4!}\sspp \frac{t\,(1\sspp 22\,t \sspp 58\,t^2\sspp 24\,t^3)}{(1\sspm t)^9}\,\frac{y^5}{5!}\sspp ...\, .
\Eeq
The coefficients of \dstyle{\frac{y^{d+1}}{(d+1)!}}, for $d\sspgeq 0$, are the \ogfs of the diagonal sequences of $J$. (In \cite{Bala} $d \sspeq n-1$.) \Eg for $d=2$, \dstyle{GDJ(2,\,t)\sspeq \frac{t\,(1\sspp 2\,t)}{(1-t)^5}} generates the third diagonal sequence $\{0,\, 1,\, 7,\, 25,\, 65,\, 140,\, 266,\, 462,\, 750,\, ...\}$ which is \seqnum{A001296}. The coefficients of the numerator polynomials are $[[1],\,[0,\,1],[0,\,1,\,2],],...$. Without the first column and offset $1$ this is \seqnum{A008517} (or \seqnum{A201637}), the second-order {\sl Eulerian} triangle, call it ${\sl Euler}2$.
\pbn  
{\bf 2)}\ $\bf P\cdot S2$:\  $S\sspeq (e^s,\, e^s\sspm 1)$. This is the product of the {\sl Sheffer} matrices $P\sspeq (e^s,\,s)$ (of the {\sl Appell} type), the {\sl Pascal} triangle \seqnum{A007318}, and $J\sspeq (1,\, e^s\sspm 1)$, {\sl Stirling}$2$ from the previous example.\pn 
Remember that {\sl Sheffer} matrices build a group (for the group law, see, \eg \cite{WLang}, {\it Lemma 9}, \Eq{139}).\pn 
Here $H(x) = \int dx e^x \sspeq e^x$, $H(0)\sspeq 1$ and the compositional inverse $x(t;\,y)$ is the one from the previous example. Now from \Eq{9}
\Beqarray
 &&EGDS(y,t)\sspeq  e^{x(t;\, y)}\sspm 1 \nonumber \\
 &\sspeq&  \frac{1}{1\sspm t}\,y \sspp \frac{1}{(1\sspm t)^3}\,\frac{y^2}{2!} \sspp \frac{1\sspp 2\,t}{(1\sspm t)^5}\,\frac{y^3}{3!} \sspp \frac{1 \sspp 8\,t \sspp 6\,t^2}{(1\sspm t)^7}\,\frac{y^4}{4!}\sspp \frac{1\sspp 22\,t \sspp 58\,t^2\sspp 24\,t^4}{(1\sspm t)^9}\,\frac{y^5}{5!}\sspp ...\, . 
\Eeqarray
This is similar to the above \egf but now the coefficient triangle for the numerator polynomials of the \ogfs is really \seqnum{A201867} (with the main diagonal $\{1,{\rm repeat}\ 0\}$).\psn 
In this way the {\sl Sheffer} triangle  $\bf PS\cdot S2$ maps to the {\sl Euler}$2$ triangle  \seqnum{A201867} (which is not {\sl Sheffer}). 
\pbn
{\bf 3)} $\bf P\cdot |S1|$:\  $S\sspeq (e^s,\, -\log(1\sspm s))$. This is the product of the {\sl Sheffer} matrices $P\sspeq (e^s,\,s)$ (of the {\sl Appell} type), the {\sl Pascal} triangle \seqnum{A007318}, and $J\sspeq (1,\, -\log(1 - s))\sspeq |{\sl Stirling}1|$ given in \seqnum{A132393}$ \sspeq |$\seqnum{A048994}$|$. This forms the {\sl Sheffer} triangle \seqnum{A094816} (coefficients of the {\sl Charlier} polynomials, see \eg \cite{Chihara}).\pn 
Here $H(x) = \int dx\, e^x \sspeq e^x$, $H(0)\sspeq 1$, like in the previous example, and the compositional inverse $x(t;\,y)$ of \dstyle{y=t(t;x)\sspeq x\,\left (1\sspm t\,\left(\frac{-\log(1\sspm x)}{x}\right)\right)} is (for Maple the expansion up to a certain power is taken)
\Beq
x(t;\,y)\sspeq \frac{1}{1\sspm t}\, y\sspp \frac{t}{(1\sspm t)^3}\, \frac{y^2}{2!}\sspp \frac{t\,(2\sspp t)}{(1\sspp t)^5}\, \frac{y^3}{3!}\sspp \frac{t\,(6\spp 8\,t\sspp t^2)}{(1\sspm t)^7}\,\frac{y^4}{4!}\sspp  \frac{t\,(24 \sspp 58\,t \sspp 22\,t^2\sspp t^3)}{(1\sspm t)^9}\, \frac{y^5}{5!}\sspp ...\, .
\Eeq
Compare this with the different \Eq{11}. Now from \Eq{9}, 
\Beqarray
EGDS(y,t)&\sspeq&  e^{x(t;\, y)}\sspm 1 \sspeq  \frac{1}{1\sspm t}\,y \sspp \frac{1}{(1\sspm t)^3}\,\,\frac{y^2}{2!} \sspp  \frac{1\sspp 3\,t\sspm t^2}{(1\sspm t)^5}\,\frac{y^3}{3!} \sspp \nonumber \\
 &\sspp& \frac{t \sspp 17\,t \sspm 2\,t^2\sspm t^3}{(1\sspm t)^7}\,\frac{y^4}{4!}\sspp \frac{1\sspp 80\,t \sspp 49\,t^2\sspm 27\,t^3\sspp 2\,t^4}{(1\sspm t)^9}\,\frac{y^5}{5!}\sspp ...\, . 
\Eeqarray
The coefficients of the row polynomials are given as signed triangle \seqnum{A290311}. Like  $\bf P\cdot S2$ produced the {\sl Euler}$2$ triangle in example $2$, here $\bf P\cdot |S1|$ produces triangle \seqnum{A290311}.
\pbn
{\bf 4)} $\bf S2[d,a]$.\  $S\sspeq (e^{a\,s},\,e^{d\,s}\sspm 1)$, generalized {\sl Stirling}$2$ number triangles (\cite{WLang2}, also with references). Here $d\sspin \mathbb N_0$, $a \sspin \mathbb N_0$ and $\gcd(d,\, a)\sspeq 1$, and for $d=1$ one puts $a\sspeq 0$. Example $1$ is the instance $[d,a]\sspeq [1,0]$, and we consider here only $d\sspgeq 2$ (\ie $a\sspneq 0$). Example $2$ would appear as $d\sspeq a\sspeq 1$.
\pn
\dstyle{y(d;t;x)\sspeq x\,\left(1\sspm t\,\frac{e^{d\,x}\sspm 1}{x}\right)} with the compositional inverse $x(d;t;y)$. \dstyle{H(a;x) = \int dx\, e^{a\,x}\sspeq \frac{1}{a}\,e^{a\,x}}, \dstyle{H(a;0)\sspeq \frac{1}{a}}. From  \Eq{9}
\Beq
EGDS2(d,a;y,t)\sspeq \frac{1}{a}\,\left(e^{a\,x(d;t;y)}\sspm 1\right)\, .
\Eeq
We consider two instances.\psn
$\boldsymbol{\alpha})$\ $S\sspeq  S2[2,1]\sspeq$\seqnum{A154537}. \psn
\Beqarray
&&EGDS2(2,1;y,t) \sspeq  e^{x(2;t;y)}\sspm 1\sspeq \frac{1}{1\sspm 2\,t}\,y \sspp \frac{1\sspp 2\,t}{(1\sspm 2\,t)^3}\,\,\frac{y^2}{2!} \sspp  \frac{1\sspp 16\,t\sspp 12\,t^2}{(1\sspm 2\,t)^5}\,\frac{y^3}{3!}   \nonumber \\
&\sspp& \frac{1 \sspp 66\,t \sspp 284\,t^2\sspp 120\,t^3}{(1\sspm 2\,t)^7}\,\frac{y^4}{4!}\sspp \frac{1\sspp 224\,t \sspp 2872\,t^2\sspp 5952\,t^3\sspp 1680\,t^4}{(1\sspm 2\,t)^9}\,\frac{y^5}{5!}\sspp ...\, .
\Eeqarray
The coefficients of the numerator polynomials are found in \seqnum{A290315}.
\psn
$\boldsymbol{\beta})$\ $S\sspeq  S2[3,1]\sspeq$\seqnum{A282629}. \psn
\Beqarray
&&EGDS2(3,1;y,t) \sspeq  e^{x(3;t;y)}\sspm 1\sspeq \frac{1}{1\sspm 3\,t}\,y \sspp \frac{1\sspp 3\,t}{(1\sspm 3\,t)^3}\,\,\frac{y^2}{2!} \sspp  \frac{1\sspp 16\,t\sspp 12\,t^2}{(1\sspm 3\,t)^5}\,\frac{y^3}{3!}   \nonumber \\
&\sspp& \frac{1 \sspp 66\,t \sspp 284\,t^2\sspp 120\,t^3}{(1\sspm 3\,t)^7}\,\frac{y^4}{4!}\sspp \frac{1\sspp 224\,t \sspp 2872\,t^2\sspp 5952\,t^3\sspp 1680\,t^4}{(1\sspm 3\,t)^9}\,\frac{y^5}{5!}\sspp ...\, .
\Eeqarray
The coefficients of the numerator polynomials are found in \seqnum{A290316}.
\pbn
{\bf 5)} $\bf \widehat{S1p}[d,a]$.\  $S\sspeq ((1\sspm d\,s)^{-\frac{a}{d}},\,-\frac{1}{d}\,\log(1\sspm d\,s)$, generalized signless {\sl Stirling}$1$ number triangles (see \cite{WLang2}, also with references). Here $d\sspin \mathbb N_0$, $a \sspin \mathbb N_0$ and $\gcd(d,\, a)\sspeq 1$, and for $d=1$ one puts $a\sspeq 0$. The $[d,a]\sspeq [1,0]$ case has been given for the signed {\sl Stirling}$1$ numbers in the Bala article \cite{Bala}, and we consider here only $d\sspgeq 2$ (\ie $a\sspneq 0$). \pn
\dstyle{y(d;t;x)\sspeq x\,\left(1\sspm t\,\left(-\frac{\log(1\sspm d\,x)}{d\,x}\right)\right)} with the compositional inverse $x(d;t;y)$. No confusion with above $y$ and $x$ quantities with the same name should arise. \pn
\dstyle{H(d,a;x) = \int dx\,(1\sspm d\,x)^{-\frac{a}{d}} \sspeq -\frac{1}{d\sspm a}\,\left(1\sspm d\,x\right)^{\frac{d-a}{d}}}, \dstyle{H(d,a;0)\sspeq -\frac{1}{d\sspm a}}. From  \Eq{9}
\Beq
EGD\widehat{S1p}(d,a;y,t)\sspeq \frac{1}{d\spm a}\,\left[1 \sspm \left(1\sspm d\,x(d;t;y)\right)^{\frac{d-a}{a}}\right]\, .
\Eeq
We consider two instances.\psn
$\boldsymbol{\alpha})$\ $S\sspeq  \widehat{S1p}[2,1]\sspeq$\seqnum{A028338}. \psn
\Beqarray
EGD\widehat{S1p}(2,1;y,t) &\sspeq& 1 - (1 - 2\,x(2;t;y))^{1/2}\sspeq \frac{1}{1\sspm t}\,y \sspp \frac{1\sspp t}{(1\sspm t)^3}\,\,\frac{y^2}{2!} \sspp  \frac{3\sspp 8\,t\sspp t^2}{(1\sspm t)^5}\,\frac{y^3}{3!}   \nonumber \\
&\sspp& \frac{15 \sspp 71\,t \sspp 33\,t^2\sspp t^3}{(1\sspm t)^7}\,\frac{y^4}{4!}\sspp \frac{105\sspp 744\,t \sspp 718\,t^2\sspp 112\,t^3\sspp t^4}{(1\sspm t)^9}\,\frac{y^5}{5!}\sspp ...\, .
\Eeqarray
The coefficients of the numerator polynomials are found in \seqnum{A288875}. The first diagonal sequences of \seqnum{A028338} are  \seqnum{A000012}, \seqnum{A000290}$(n+1)$, \seqnum{A024196}$(n+1)$, \seqnum{A024197}$(n+1)$, \seqnum{A024198}$(n+1)$.
\pbn
$\boldsymbol{\beta})$\ $S\sspeq  \widehat{S1p}[3,1]\sspeq$\seqnum{A286718}. \psn
\Beqarray
EGD\widehat{S1p}(3,1;y,t)\sspeq (1 - (1 - 3\,x(3;t;y))^{2/3})/2\sspeq \frac{1}{1\sspm t}\,y \sspp \frac{1\sspp 2\,t}{(1\sspm t)^3}\,\,\frac{y^2}{2!} \sspp  \frac{4\sspp 19\,t\sspp 4\,t^2}{(1\sspm t)^5}\,\frac{y^3}{3!}   &&\nonumber \\
\sspp \frac{28 \sspp 222\,t \sspp 147\,t^2\sspp 8\,t^3}{(1\sspm t)^7}\,\frac{y^4}{4!}\sspp \frac{280\sspp 3194\,t \sspp 4128\,t^2\sspp 887\,t^3\sspp 16\,t^4}{(1\sspm t)^9}\,\frac{y^5}{5!}\sspp ...\, .&&
\Eeqarray
The coefficients of the numerator polynomials are found in \seqnum{A290318}. The first diagonal sequences of \seqnum{A286718} are \seqnum{A000012}, \seqnum{A000326}$(n+1)$, \seqnum{A024212}$(n+1)$, \seqnum{A024213}$(n+1)$.
\pbn
\pbn
{\bf Part B: Riordan triangles and their diagonals multiplied with Pascal diagonals} 
\pbn 
\hskip 1cm A {\sl Riordan} triangle $R$ (an infinite dimensional lower triangular (ordinary) convolution matrix; for practical purpose a $N\times N$ matrix) is denoted by $R\sspeq (G,\, F)$ with \ogf \dstyle{G(x)\sspeq \sum_{k=0}^{\infty}\, G_n\,x^n}, where $G(0) \sspeq G_0 \speq 1$ ({\it w.l.o.g.}), and  $F(x) \sspeq x\,{\widehat F}(x)$ with \ogf \dstyle{{\widehat F}(x)\sspeq \sum_{k=0}^{\infty}\, {\widehat F}_n\,x^n}, where ${\widehat F}(0)\sspeq {\widehat F}_0\sspneq 0$. The column sequence $RCol(m)\sspeq \{R(n, m)\}_{n=0}^{\infty}$ (with $m$ leading zeros) has \ogf \dstyle{GRCol(x, m)\sspeq \sum_{n = m}^{\infty}\, R(n, m)\,x^n}, for $m\sspin \mathbb N_0$, given by
\Beq
GRCol(x, m)\sspeq G(x)\,F(x)^m\sspeq G(x)\,x^m\,{\widehat F}(x)^m\, .
\Eeq
The row polynomials (called {\sl Riordan} polynomials) are  $PR(n,\, x)\sspeq \sum_{m=0}^n\,R(n,\,m)\,x^n$. They have \ogfs \dstyle{GPR(x,z)\sspeq \sum_{n=0}^{\infty} PR(n,\, x)\,z^n} given by
\Beq
GPS(x,z)\sspeq G(z)\,\frac{1}{1\sspm x\,F(z)}\, 
\Eeq
which is also called the \ogf of the $R$ triangle. \psn
The {\sl Riordan} group has been introduced, in analogy to the {\sl Sheffer} group \cite{Roman} by {\sl Shapiro} {\it et al.} \cite{Shapiro}\psn 
There is no (ordinary) convolution property for {\sl Riordan} polynomials similar to \Eq{3}. But $P\sspeq (1,\, F)$ is also called associated to $R\sspeq (g,\, f)$. Such matrices form a subgroup of the {\sl Riordan} group.\psn
The diagonal sequences are labeled by $d\sspin \mathbb N_0$, with $d\sspeq 0$ for the main diagonal. Their entries are
\Beq
DR(d,\,m)\sspeq R(d+m,\,m),\ \ {\rm for}\ \ m\sspin \mathbb N_0\,.
\Eeq
Their \ogf is
\Beq
GDR(d,x)\sspeq \sum_{m=0}^{\infty}\, DR(d,\,m)\,x^m\,,
\Eeq
Application of {\sl Lagrange}'s theorem, like in the {\it Lemma}, part b) does not lead to the \ogfs of these diagonal sequences directly. Instead one is led to consider the product of the diagonal entries with the corresponding ones of {\sl Pascal}'s {\sl Riordan} triangle \dstyle{P\sspeq \left(\frac{1}{1\sspm x},\,\frac{x}{1\sspm x}\right)}, \seqnum{A007318}. This belongs to the so called {\sl Bell} subgroup of the {\sl Riordan} group of the type $B\sspeq (G(x),\, x\,G(x))$.  Define
\Beq
{\widehat D}(d,\,m)\spdef P(d+m, m)\,R(d+m,\,m)\sspeq {\binomial{d+m}{m}}\,D(d,\,m)\, .
\Eeq
The corresponding \ogf is $G{\widehat D}(d,\,t)\sspeq \sum_{m=0}^{\infty}\,{\widehat D}(d,\,m)\,t^m $
Their logarithmic generating function (\lgf) $LG{\widehat D}R(y,t)$ is taken as
\Beq
LG{\widehat D}R(y,\,t)\sspeq \sum_{d=0}^{\infty}\,G{\widehat D}(d,\,t)\,\frac{y^{d+1}}{d+1}\, .
\Eeq
(The unconventional powers for this \lgf and the use of $y$ instead of $z$ will become clear later). \psn
Applying now {\it Lemma}, part b) to $y\sspeq y(t;x)\sspeq x\,\psi(t;x)\sspeq x\,(1\sspm t\, {\widehat F}(x))\sspeq x\sspm t\, F(x)$ with the {\sl Riordan} function $F$, introducing a parameter $t$, and taking $H(x)\sspeq \int\,dx\,G(x)$, with the {\sl Riordan} function $G$, we obtain, with the compositional inverse $x\sspeq x(t;y)$ of $y\sspeq y(t;x)$, the following proposition. 
\pbn
{\bf Proposition 2:}
\Beq
LG{\widehat D}R(y,t) \sspeq H(x(t;y))\sspm H(0)\sspeq \left.\left[\int\,dx\,G(x)\right]\right|_{x\sspeq x(t;y)}\sspm \left.\left [\int\,dx\,G(x)\right]\right|_{x\sspeq 0}\ .
\Eeq
As in the {\sl Sheffer} section one has first to compute the  $x\sspeq x(t;y)$, the {\sl Lagrange} inversion of $y\sspeq y(t;x)$. This is the case $H(x) =x$ in the {\it Lemma}, part b, with the chosen $\psi\sspeq \psi(t;x)\sspeq 1-t\,{\widehat F}(x)$. It belongs to the associated {\sl Riordan} case $A\sspeq (1,\,F)$ ({$A$ for the associated triangle to $R$). This yields the following corollary.\psn
{\bf Corollary 2: Associated Riordan case}\psn
\Beq
LG{\widehat D}A(y,t) \sspeq x(t;y)\, .
\Eeq
This means that in the $A\sspeq (1,\, F)$ case the \lgf of the \ogfs of the sequences of the product of the entries of the diagonals of $A$ and the {\sl Pascal} triangle $P$ is just the compositional inverse of $y\sspeq y(t;x)\sspeq x\sspm t\,F(x)$.\psn
Instead of the \lgf of the \ogfs of diagonal sequences of the triangle with entries  ${\widehat D}(d,\,m)$ one could as well take the \egf  of the \egfs of the diagonal sequences of the triangle with entries ${\widetilde D}(d,\,m)\sspdef (d+m)!\,D(d,\,m)$. This leads to
\psn
{\bf Corollary 3:} \ With the \egf
\Beq
E{\widetilde D}(d,\,t)\sspdef \sum_{m=0}^{\infty}\, {\widetilde D}(d,\,m)\,\frac{t^m}{m!}\sspeq \sum_{m=0}^{\infty}\,(d+m)!\,D(d,\,m)\,\frac{t^m}{m!}\, ,
\Eeq
and the further \egf
\Beq
EE{\widetilde D}R(y,\,t)\sspeq \sum_{d=0}^{\infty}\,E{\widetilde D}(d,\,t)\,\frac{y^{d+1}}{(d+1)!}
\Eeq
one has
\Beq
EE{\widetilde D}R(y,\,t)\sspeq H(x(t;y))\sspm H(0)\sspeq \left.\left[\int\,dx\,G(x)\right]\right|_{x\sspeq x(t;y)}\sspm \left.\left [\int\,dx\,G(x)\right]\right|_{x\sspeq 0}\ .
\Eeq
\pbn
{\bf Examples}
\pbn
{\bf 1)} \dstyle{A\sspeq \left(1,\,\frac{x}{1\sspm x}\right)}, the {\bf Pascal} triangle variant given in \seqnum{A097805}.
The {\sl Lagrange} inverse  $x\sspeq x(t;y)$ of \dstyle{y\sspeq x\,\left (1\sspm \frac{t}{1\sspm x}\right)} turns out to be 
\Beq
 x(t;y)\sspeq \frac{1}{1\sspm t}\,y \sspp \frac{2\,t}{(1\sspm t)^3}\, \frac{y^2}{2} \sspp \frac{3\,t\,(1\sspp t)}{(1\sspm t)^5}\,\frac{y^3}{3}\sspp \frac{4\,t\,(1 \sspp 3\,t \sspp t^2)}{(1\sspm t)^7}\,\frac{y^4}{4}\sspp \frac{5\,t\,(1\sspp 6\,t \sspp 6\,t^2\sspp t^3)}{(1\sspm t)^9}\,\frac{y^5}{5}\sspp ...\, .
\Eeq
See \cite{Drake}, Example $1.10.8$.\pn
This is a \lgf, therefore the coefficients of \dstyle{\frac{y^{d+1}}{d+1}}, for $d\sspgeq 0$, are the \ogf of the diagonal sequences of the triangle $[[1], [0, 1], [0, 2, 1], [0, 3, 6, 1], [0, 4, 18, 12, 1], [0, 5, 40, 60, 20, 1], [0, 6, 75, 200, 150, 30, 1], [0, 7, 126, 525, 700, 315, 42, 1], ...]$ obtained by multiplying the entries of {\sl Pascal}'s triangle and $A\sspeq$\seqnum{A097805}. \Eg the fourth diagonal ($d=3$) $[0,\, 4,\,40,\, ...]$ has o.g.f. \dstyle{G(3, x)  =  \frac{4\,t\,(1 \sspp 3\,t \sspp t^2)}{(1\sspm t)^7}}. The numerator polynomials divided by $(d+1)\,t$, for $d\sspgeq 1$, are found as row $d$ polynomials of \seqnum{A001263} ({\sl Narayana} triangle).
\pbn
{\bf 2)\ Generalized Pascal triangles.}\psn
\dstyle{R\sspeq \left(G(x),\, \frac{x}{1\sspm x}\right)}, and the {\sl Lagrange} inverse $x(t;y)$ is given in \Eq{32}.
Now \Eq{27} applies with $H(x)\sspeq \int dx\,G(x)$.\psn
Two instances:\psn
$\boldsymbol{\alpha})$\ \dstyle{R\sspeq \left(\frac{1}{1\sspm x},\,\frac{x}{1\sspm x}\right)}. This is the {\sl Pascal} triangle \seqnum{A007318}. Here $H(x)\sspeq -\log(1\sspm x)$, $H(0)\sspeq 0$, and one obtains the \lgf
\Beqarray
LG{\widehat D}R(y,t) &\sspeq& -\log(1\sspm x(y;t))\sspeq  \frac{1}{1\sspm t}\,y \sspp \frac{1\sspp t}{(1\sspm t)^3}\, \frac{y^2}{2} \sspp \frac{1\sspp 4\,t\sspp t^2}{(1\sspm t)^5}\,\frac{y^3}{3}\sspp \nonumber \\
&\sspp&  \frac{1 \sspp 9\,t \sspp 9\,t^2\sspp t^3}{(1\sspm t)^7}\,\frac{y^4}{4}\sspp \frac{1\sspp 16\,t \sspp 36\,t^2\sspp 16\,t^3\sspp t^4}{(1\sspm t)^9}\,\frac{y^5}{5}\sspp ...
\Eeqarray
The numerator polynomials are the row polynomials of \seqnum{A008459}, the square entries of {\sl Pascal}'s triangle.
The \ogfs for the diagonal sequences  of \seqnum{A008459} are given by \dstyle{GDR(d,\,x)\sspeq \left[\frac{y^{d+1}}{d+1}\right]\, LG{\widehat D}R(y,t)} for $d\sspgeq 0$. \Eg the fourth diagonal sequence $[1,\,16,\,100,\,...]$ has \ogf \dstyle{GDR(3,\,x)\sspeq \frac{1 \sspp 9\,t \sspp 9\,t^2\sspp t^3}{(1\sspm x)^7}}.
\pbn 
$\boldsymbol{\beta})$\  \dstyle{R\sspeq \left(\frac{1}{(1\sspm x)^2},\,\frac{x}{1\sspm x}\right)}. This is the {\sl Riordan} triangle \seqnum{A135278}.
Here \dstyle{H(x)\sspeq \frac{1}{1\sspm x}}, $H(0)\sspeq 1$, and one obtains the \lgf
\Beqarray
LG{\widehat D}R(y,t) &\sspeq& \frac{1}{1\sspm x(y;t)}\sspm 1\sspeq  \frac{1}{1\sspm t}\,y \sspp \frac{2}{(1\sspm t)^3}\, \frac{y^2}{2} \sspp \frac{3\,(1\sspp t)}{(1\sspm t)^5}\,\frac{y^3}{3}\sspp \nonumber \\
&\sspp&  \frac{4\,(1 \sspp 3\,t \sspp t^2)}{(1\sspm t)^7}\,\frac{y^4}{4}\sspp \frac{5\,(1\sspp 6\,t \sspp 6\,t^2\sspp t^3)}{(1\sspm t)^9}\,\frac{y^5}{5}\sspp ...
\Eeqarray
The numerator polynomials are again the row polynomials of \seqnum{A008459} ({\sl Narayana} triangle) multiplied here by $d+1$.
Therefore, the \ogfs for the diagonal sequences with entries \seqnum{A103371}$(n,\,k)\sspeq$ \seqnum{A135278}$(n,\,k)$\,\seqnum{A007318}$(n,\,k)$  are given by \dstyle{GDR(d,\,x)\sspeq (d+1)\,\frac{\sum_{k=1}^d\,N(d,\,k)\,x^{k-1}}{(1\sspm x)^{2\,d+1}}} for $d\sspgeq 1$, with $N(d,\,k)\sspeq$\seqnum{A008459}$(d,\, k)$, and for $d\sspeq 0$ the \ogf is \dstyle{GDR(0,\,x)\sspeq \frac{1}{1\sspm x}}. 
\pbn
\pbn
\section{Proofs}
{\bf Part A}\psn
{\bf 1.}\  {\bf Proof of the Lemma: Lagrange theorem and inversion} \cite{Fichtenholz}, p. 523. \Eq{29}, \cite{WhittakerWatson}, p. 133.\psn
Part {\bf a)} is the standard theorem of {\sl Lagrange} with the proof given in the references.\psn
Part {\bf b)}:\ The first two equations of \Eq{8} follow from part {\bf a)} for $a\sspeq 0$, interchanging the r\^ole of $x$ and $y$, and using \dstyle{\varphi(x)\sspeq \frac{1}{\psi(x)}} (See\cite{Fichtenholz}, pp. 524-525 for the case $H(x)\sspeq x$). The last eq. is then obvious with the definition of $[a^n]\,h(a)$ given there.
\psn
{\bf 2.}\  {\bf Proof of Proposition 1}\psn
From the {\it Lemma}, part b), one has, with $y\sspeq y(t;x)\sspeq x\, \psi(t;x)\sspeq x\,(1\sspm t\,{\hat f}(x))$, and $H(x)\sspeq \int dx\, g(x)$, where the {\sl Sheffer} triangle is $S\sspeq (g(x),\, x\,{\hat f}(x))$,
\Beq
H(x(t;y))\sspm H(0)\sspeq \sum_{n=1}^{\infty}\, \frac{y^n}{n!}\,(n-1)!\,[a^{n-1}]\,\left[(1\sspm t\,{\hat f}(a))^{-n}\,g(a)\right]\, .
\Eeq 
The binomial theorem \dstyle{(1\sspm t\,{\hat f}(a))^{-n}\sspeq \sum_{p=0}^{\infty}\, {\binomial{-n}{p}}\, (-t)^p\,({\hat f}(a))^p  } is applied. Then the binomial with negative upper entry is transformed in one with non-negative entries, using the identity
(see \cite{GKP}, p. 164, \Eq{5.14})
\Beq
{\binomial{-n}{p}}\sspeq (-1)^p\,{\binomial{p+n-1}{p}}\ .
\Eeq
\Beq
H(x(t;y))\sspm H(0)\sspeq \sum_{n=1}^{\infty}\, \frac{y^n}{n!}\,(n-1)!\,\sum_{p=0}^{\infty}\, {\binomial{p+n-1}{p}}\,t^p\,p!\,[a^{n-1}]\,\left[\frac{({\hat f}(a))^p}{p!}\,g(a)\right]\, .
\Eeq 
In order to obtain $f(a)\sspeq a\,{\hat f}(a)$ one uses $[a^{n-1}]\,h(a)\sspeq [a^{n-1+p}](a^p\,h(a))$. Then the definition of the \egf of the sequence of column $p$ of the {\sl Sheffer} triangle is used: \dstyle{\sum_{k=p\, (0)}^{\infty} S(k,\,p)\,\frac{a^k}{k!}\sspeq \frac{(f(a))^p}{p!}\, g(a)} (One can start with $p=0$ because $S(k,\,p)\sspeq 0$ for $0\sspleq k\sspkl p$.) Thus \dstyle{[a^m] \left(\frac{(f(a))^p}{p!}\, g(a)\right) \sspeq S(m,\, p)\,\frac{1}{m!}}.
\Beqarray
H(x(t;y))\sspm H(0)&\sspeq& \sum_{n=1}^{\infty}\, \frac{y^n}{n!}\,(n-1)!\, \sum_{p=0}^{\infty}\, {\binomial{p+n-1}{p}}\,t^p\,p!\, \frac{1}{(n-1+p)!}\,S(n-1+p,\,p) \nonumber \\
&\sspeq& \sum_{n=1}^{\infty}\, \frac{y^n}{n!}\,\left(\sum_{p=0}^{\infty}\, t^p\,S(n-1+p,\,p)\right)\, .
\Eeqarray
But the \ogf of the diagonal sequences of $S$ is $GDS(n-1,\,t)\sspeq \sum_{p=0}^{\infty}\, t^p\,S(n-1+p,\,p)$, for $n\sspgeq 1$, and because we take $d\sspeq n-1$ to label the diagonals, we get 
\Beq
H(x(t;y))\sspm H(0) \sspeq \sum_{d=0}^{\infty}\, \frac{y^{d+1}}{(d+1)!}\, GDS(d,\,t)\sspfed EGDS(y,\,t)
\Eeq
\hskip 17.1cm $\square$
\pbn
{\bf Part B}
{\bf 3.} {\bf Proof of Proposition 2}\psn
From the {\it Lemma}, part b), one has, with $y\sspeq y(t;x)\sspeq x\, \psi(t;x)\sspeq x\,(1\sspm t\,{\widehat F}(x))$, and $H(x)\sspeq \int dx\, G(x)$, where the {\sl Riordan} triangle is $R\sspeq (G(x),\, x\,{\widehat F}(x))$
\Beq
H(x(t;y))\sspm H(0)\sspeq \sum_{n=1}^{\infty}\, \frac{y^n}{n!}\,(n-1)!\,[a^{n-1}]\,\left[(1\sspm t\,{\widehat F}(a))^{-n}\,G(a)\right]\, .
\Eeq 
Using the binomial theorem and the binomial identity \Eq{36} one finds 
\Beq
H(x(t;y))\sspm H(0)\sspeq \sum_{n=1}^{\infty}\, \frac{y^n}{n!}\,(n-1)!\,\sum_{p=0}^{\infty}\, {\binomial{p+n-1}{p}}\,t^p\,[a^{n-1}]\,\left[({\widehat F}(a))^p\,G(a)\right]\, .
\Eeq 
In order to obtain $F(a)\sspeq a\,{\widehat F}(a)$ one uses $[a^{n-1}]\,h(a)\sspeq [a^{n-1+p}](a^p\,h(a))$. Then the definition of the \ogf of the sequence of column labeled $p$ of triangle $R$ is used: \dstyle{\sum_{k=p\, (0)}^{\infty} R(k,\,p)\,a^k\sspeq (F(a))^p\, G(a)}. Thus \dstyle{[a^m] \left((F(a))^p\, G(a)\right) \sspeq R(m,\, p)}.
\Beqarray
H(x(t;y))\sspm H(0)&\sspeq& \sum_{n=1}^{\infty}\, \frac{y^n}{n!}\,(n-1)!\, \sum_{p=0}^{\infty}\, {\binomial{p+n-1}{p}}\,t^p\, R(n-1+p,\,p) \nonumber \\
&\sspeq& \sum_{n=1}^{\infty}\, \frac{y^n}{n!}\,\sum_{p=0}^{\infty}\, \frac{t^p}{p!}\,(p+n-1)!\,\,R(n-1+p,\,p)\, .
\Eeqarray
At this stage the {\it Corollary 3} has been proved, if one uses $n-1\sspeq d$ (and $p\sspto m$). But we prefer to use the binomial coefficient to multiply the diagonal $R$ entries, \ie we use the first equation. With $n-1\sspeq d$ this becomes
\Beq
H(x(t;y))\sspm H(0)\sspeq \sum_{d=0}^{\infty}\, \frac{y^{d+1}}{d+1}\,\left(\sum_{p=0}^{\infty}\,{\binomial{d+p}{p}}\,R(d+p,\,p)\,t^p \right)\,.
\Eeq
This is the the \lgf \Eq{26} of the \ogfs  $G{\widehat D}(d, t)$ of the product of the diagonal entries in {\sl Pascal}'s triangle and the ones of the {\sl Riordan} triangle, called ${\widehat D}(d\, m)$ in \Eq{25}. \hskip 9cm $\square$
\pbn
\vfill\eject\noin

\pbn
\hrulefill
\pbn
{\it 2010 Mathematics Subject Classification}: Primary 05A15, 11B83, Secondary 11B37.\psn
{Keywords}: Generating functions (ordinary, exponential, logarithmic), Sheffer arrays, Riordan arrays, Psacal triangle, Stirling triangles.
\pbn
\hrulefill
\pbn
{\it OEIS} \cite{OEIS} A numbers: \psn
\seqnum{A000012}, \seqnum{A000290}, \seqnum{A000326}, \seqnum{A001263}, \seqnum{A001296}, \seqnum{A006232}, \seqnum{A007318}, \seqnum{A008459}, \seqnum{A008517},  \seqnum{A024196}, \seqnum{A024197}, \seqnum{A024198}, \seqnum{A024212}, \seqnum{A024213}, \seqnum{A028338}, \seqnum{A048993}, \seqnum{A048994},  \seqnum{A094816}, \seqnum{A097805}, \seqnum{A103371}, \seqnum{A112007}, \seqnum{A132393}, \seqnum{A135278}, \seqnum{A154537}, \seqnum{A201637}, \seqnum{A201867}, \seqnum{A282629}, \seqnum{A286718}, \seqnum{A288875}, \seqnum{A290311}, \seqnum{A290315}, \seqnum{A290316}, \seqnum{A290318}.
\pbn
\hrulefill
\end{document}